\documentclass[reqno,a4paper]{amsart}

\usepackage{amsmath, amsfonts, amsthm, amssymb,stmaryrd}
\usepackage{array}
\usepackage{longtable} 
\usepackage{colortab} 
\usepackage{colortbl}
\usepackage{arydshln}
\usepackage{dsfont}
\usepackage{graphicx} 
\usepackage{boxedminipage}
\usepackage[left=2cm, right=2cm, top=3cm,bottom=3cm]{geometry}
\usepackage[all]{xy}
\usepackage{wrapfig}
\usepackage{enumerate}
\usepackage{amssymb}
\usepackage{graphicx,color} 
\usepackage{tikz}
\usepackage[english]{babel}
\usepackage[utf8]{inputenc}
\usepackage[T1]{fontenc}
\usepackage{ upgreek }
\usepackage{xcolor}
\usepackage{setspace}
\usepackage{url}
\usepackage{hyperref}

\theoremstyle{definition}
\newtheorem{teo}{Theorem}[section]
\newtheorem{prop}[teo]{Proposition}
\newtheorem{lem}[teo]{Lemma}
\newtheorem{conj}[teo]{Conjecture}

\def\d{{\sf d}}
\def\D{{\sf D}}

\def\v{{\sf v}}
\def\L{{\sf L}}
\def\HH{{\mathbb M}}
\def\LL{{\mathcal L}}
\def\F{{\mathcal F}}
\def\B{{\mathcal B}}
\def\U{{\mathcal U}}
\def\A{{\mathcal A}}
\def\bd{{\boldsymbol{\cdot}}}
\def\supp{\text{supp}}

\def\G{\mathcal G}
\def\N{\mathbb N}

\def\Z{\mathbb Z}

\def\Q {\mathbb Q}

\def\ord{\mathop{\rm ord}\nolimits}

\newcommand{\proofthmmain}{\noindent {\it Proof of Theorem \ref{main}. }}

\title{On minimal product-one sequences of maximal length over the non-abelian group of order $pq$}

\author[D.V. Avelar]{D.V. Avelar$^1$}
\address{Departamento de Análise\\
	Universidade Federal Fluminense (UFF)\\
	Niterói, RJ\\
	24210-201\\
	Brazil\\
}
\email{daniloavelar@id.uff.br}

\author[F.E. Brochero Mart\'{\i}nez]{F.E. Brochero Mart\'{\i}nez$^2$}
\address{
	Departamento de Matem\'{a}tica\\
	Universidade Federal de Minas Gerais (UFMG)\\
	Belo Horizonte, MG\\
	31270-901\\
	Brazil\\
}
\email{fbrocher@mat.ufmg.br}

\author[S. Ribas]{S. Ribas$^3$}
\address{
	Departamento de Matem\'{a}tica\\
	Universidade Federal de Ouro Preto\\
	Ouro Preto, MG\\
	35402-136\\
	Brazil
}
\email{savio.ribas@ufop.edu.br}

\thanks{$^1$partially supported by FAPERJ grant E-26/210.315/2024 and FAPEMIG grant APQ-04712-25. \\
	$^2$partially supported by FAPEMIG grant RED-00133-21, and CNPq grants 316843/2023-7 and 420721/2025-8. \\	
	$^3$partially supported by FAPEMIG grants RED-00133-21, APQ-01712-23 and APQ-04712-25, and CNPq grant 420721/2025-8.
	}

\date{\today}

\begin{document}
	
\maketitle
	
\begin{abstract}		
	Let $G$ be a finite group. A sequence over $G$ is a finite multiset of elements of $G$, and it is called product-one if its terms can be ordered so that their product is the identity of $G$. The large Davenport constant $\D(G)$ is the maximal length of a minimal product-one sequence, that is, a product-one sequence that cannot be partitioned into two nontrivial product-one subsequences. Let $p,q$ be odd prime numbers with $p \mid q-1$ and let $C_q \rtimes C_p$ denote the non-abelian group of order $pq$. It is known that $\D(C_q \rtimes C_p) = 2q$. In this paper, we describe all minimal product-one sequences of length $2q$ over $C_q \rtimes C_p$. As an application, we further investigate the $k$-th elasticity (and, consequently, the union of sets containing $k$) of the monoid of product-one sequences over these groups.
\end{abstract}

\section{Introduction}

Let $G$ be a finite group written multiplicatively. A sequence $S$ over $G$ is a finite multiset of elements from $G$. The {\em zero-sum problems} investigate the conditions under which a given sequence over $G$ contains a subsequence whose product is the identity of $G$. Such subsequences are called {\em product-one}. The {\em small Davenport constant} $\d(G)$ is defined as the maximal length of a sequence over $G$ that do not have product-one subsequences. Moreover, the {\em large Davenport constant} $\D(G)$ is defined as the maximal length of a minimal product-one sequence, that is, a product-one sequence that cannot be decomposed into two nontrivial product-one subsequences. This invariant appears in the pioneering works due to Rogers \cite{Rog}, van Emde Boas and Kruyswijk \cite{vEBK} and Olson \cite{Ols1,Ols2}. In particular, Rogers showed that $\D(G)$ represents how many prime ideal factors a prime element can have in an algebraic number field with class group isomorphic to $G$. This is a crucial relation between zero-sum problems and factorization theory in Krull monoids (see \cite{GeRu}). 
	
By definition, it follows that $\d(G) + 1 \le \D(G)$, with equality for abelian groups. Furthermore, a simple application of Pigeonhole Principle yields $\D(G) \le |G|$, with equality for cyclic groups. Let $C_n$ be the cyclic group of order $n$. By the Fundamental Theorem of Abelian Groups, if $G$ is a nontrivial finite abelian group, then there exist unique $1 < n_1 \mid n_2 \mid \ldots \mid n_r$ such that $G \cong C_{n_1} \otimes \ldots \otimes C_{n_r}$, where $r$ is the rank of $G$ and $n_r$ is the exponent of $G$. Set $\D^*(G) = 1 + \sum_{i=1}^r (n_i-1)$. Some routine arguments yield $\D(G) \ge \D^*(G)$, with equality for abelian groups of rank at most $2$ and for $p$-groups \cite{Ols1,Ols2}. However, there exist infinitely many groups for which the inequality is strict \cite{GeSc}.
	
Although the main focus of the zero-sum problems has been on abelian groups due to the relation with factorization theory, it has been extended to non-abelian groups in the 80s. This explains the use multiplicative notation and the term {\em product-one} instead of additive notation and the term {\em zero-sum}. Nevertheless, for non-abelian groups, the small Davenport constant does not admit an interpretation in terms of monoid factorizations.
	

Concerning the large Davenport constant over non-abelian groups, we refer the reader to \cite{GeGr,Gr}. In particular, Geroldinger and Grynkiewicz \cite{GeGr} fully determined the large Davenport constant for groups having a cyclic subgroup of index $2$, and Grynkiewicz \cite{Gr} proved that if $G \cong C_q \rtimes C_p$ is the non-abelian group of order $pq$, where $p,q$ are odd primes with $p \mid q-1$, then $\D(G) = 2q$. He also obtained several upper bounds for $\D(G)$ for general groups $G$. Among others, we highlight $\D(G) \le \d(G) + 2|G'|-1$, where $G'$ is the commutator subgroup of $G$, and $\D(G) \le 2|G|/p$, where $G$ is non-cyclic and $p$ is the smallest prime divisor of $|G|$.

The {\em direct problem} associated to $\D(G)$ asks for the exact value of $\D(G)$, while the associated {\em inverse problem} seeks to describe the exceptional sequences of length exactly $\D(G)$ that cannot be partitioned into two nontrivial product-one subsequences (see also \cite{AMR1,AMR,Bas,MR1,MR2,MR3,HZ,BRGS,OhZh,OhZh1,QL,QL1,Ri,YZF} for other recent developments on the direct and inverse problems over non-abelian groups).
	
The goal of this paper is to establish, by arguments similar to those in \cite{Gr}, the inverse problem associated to the large Davenport constant of $C_q \rtimes C_p$, that is, a complete characterization of minimal product-one sequences of length $2q$ over $C_q \rtimes C_p$. This is in the same idea as \cite{OhZh1}, where Oh and Zhong solved the inverse problem for $\D(G)$ over dihedral and dicyclic groups. 
	
The minimal product-one sequences over $G$ are precisely the atoms, that is, irreducible elements of the monoid $\B(G)$ of the product-one sequences over $G$. Building on our characterization of minimal product-one sequences of length $\D(G)$, we investigate the unions of sets of lengths in $\B(G)$ by studying the $k$-th elasticity of $\B(G)$. This approach is also similar to that employed by Oh and Zhong in \cite{OhZh1} for dihedral and dicyclic groups.
	
The main result of this paper is the following.
	
\begin{teo}\label{main}
	Let $S$ be a minimal product-one sequence over $C_q \rtimes C_p$ of length $|S| = \D(C_q \rtimes C_p) = 2q$. Then there exist $x,y \in C_q \rtimes C_p$ for which $C_q \rtimes C_p = \langle x,y \colon x^p = y^q = 1, yx = xy^s, \ord_q(s)=p \rangle$ and 
	\begin{equation}\label{formaS}
		S = y^{[q-1]} \bd x \bd y^{[q-1]} \bd x^{p-1}y^{s^{p-1}+1}.
	\end{equation}
\end{teo}

The paper is organized as follows. In Section \ref{secdefs}, we present the prerequisite notation and definitions that will be used througout the paper, as well as some properties of the group $C_q \rtimes C_p$. In Section \ref{seclemas}, we present several auxiliary results required for the proof of Theorem \ref{main}. In Section \ref{secprovamain}, we prove two results that lead to the proof of Theorem \ref{main}. In Section \ref{secaplica}, we apply our description of the atoms of $C_q \rtimes C_p$ to the study of the $k$-th elasticity, which in turn yields information on the union of sets of lengths containing $k$.

\section{Notation and preliminaries}\label{secdefs}

We use the standard notation from group theory. In particular, for a finite group $G$,

\begin{itemize}
	\item if $A, B \subset G$, then the {\em product-set} of $A$ and $B$ is the set $AB = \{ab \colon a \in A, b \in B\}$. For a singleton $A = \{a\}$, we denote $aB = \{ab \colon b \in B\}$;
	\item if $A \subset G$, then $\langle A \rangle \le G$ denotes the subgroup generated by $A$;
	\item $Z(G) = \{g \in G \colon gh = hg \text{ for every } h \in G\} \trianglelefteq G$ is the {\em centre} of $G$;
	\item $[g,h] = g^{-1}h^{-1}gh \in G$ is the {\em commutator} of $g,h \in G$;
	\item $G' = \langle [g,h] \colon g,h \in G \rangle \trianglelefteq G$ is the {\em commutator subgroup} of $G$;
	\item $C_G(g) = \{h \in G \colon gh = hg\} \le G$ is the {\em centralizer} of $g \in G$; and
	\item for $A, B \subset G$ and $g \in G$, set the conjugations $A^g = \{g^{-1}ag \colon a \in A\}$ and $A^B = \{b^{-1}ab \colon a \in A, b \in B\}$.
\end{itemize}

In what follows, we present the necessary definitions concerning sequences, ordered sequences, and the group under consideration in this paper. The notation in this paper is consistent with \cite{Gr,OhZh1}.

\subsection{Sequences over groups}

Let $G$ be a finite group written multiplicatively and let $\F(G)$ be the free abelian monoid with basis $G$ with operation denoted by the bold dot $\bd$. A sequence $S$ over $G$ is an element of $\F(G)$, meaning that $S$ is a finite multiset of elements of $G$ (allowing repetition and the order is disregarded). In particular, if $S \in \F(G)$, then 
$$S = g_1 \bd \dots \bd g_k = \prod_{i=1}^k g_i = \prod_{g \in G} g^{[\v_g(S)]},$$
where $\v_g(S)$ is the {\em multiplicity} of $g$ in $S$ and $|S| = k = \sum_{g \in G} \v_g(S)$ is the {\em length} of $S$. We observe that $g_1 \cdot g_2 = g_1g_2 \in G$ denotes the product of $g_1$ and $g_2$, while $g_1 \bd g_2 \in \F(G)$ denotes a two-term sequence. 
A {\em subsequence} of $S$ is a divisor $T \mid S$ in $\F(G)$. In other words, $T \mid S$ if and only if $\v_g(T) \le \v_g(S)$ for every $g \in G$. In this case, we write $S \bd T^{[-1]} = \prod_{g \in G} g^{[\v_g(S) - \v_g(T)]}$. For a subset $K \subset G$, we denote $\v_K(S) = \sum_{g \in K} \v_g(S)$. Moreover, the {\em support} of $S$ is the set $\supp(S) = \{g \in G \colon \v_g(S) > 0\}$.

The {\em set of products} and the {\em set of subproducts} of $S$ are 
$$\pi(S) = \left\{ \prod_{i=1}^{|S|} g_{\sigma(i)} \in G \colon \sigma \text{ is a permutation of }[1,|S|]\right\} \quad \text{ and } \quad \Pi(S) = \bigcup_{T \mid S \atop |T| \ge 1} \pi(T),$$ 
respectively. The sequence $S \in \F(G)$ is called 
\begin{enumerate}[(i)]
	\item {\em trivial} if $|S| = 0$ (in this case, $S$ is the identity of $\F(G)$);
	\item {\em product-one} if $1 \in \pi(S)$;
	\item {\em product-one free} if $1 \not\in \Pi(S)$; 
	and
	\item {\em minimal product-one} if $1 \in \pi(S)$ and $S \neq T_1 \bd T_2$ for every $T_1, T_2$ nontrivial product-one sequences.
\end{enumerate} 

Let 
$$\B(G) = \{S \in \F(G) \colon 1 \in \pi(S)\}$$ 
denote the {\em set of product-one sequences over $G$}, and let $$\A(G) = \{S \in \B(G) \colon S \text{ is minimal product-one}\}.$$ 
We observe that $\B(G)$ is a submonoid of $\F(G)$ and $\A(G)$ is the set of {\em atoms} (or irreducible elements) of $\B(G)$. With this notation, the {\em large Davenport constant} is 
$$\D(G) = \sup\{|S| \colon S \in \A(G)\},$$ 
and the {\em small Davenport constant} is 
$$\d(G) = \sup\{|S| \colon S \in \F(G) \text{ and } 1 \not\in \Pi(S)\}.$$ 

Moreover, we observe that if $S \in \F(G)$, then $\pi(S)$ is contained in a $G'$-coset, that is, $\pi(S) = Ag$ for some $A \subset G'$ and some $g \in G$. This implies that if $S_1, \dots, S_t \in \F(G)$, then, for each $j \in [1,t]$, $\pi(S_j) = A_jg_j$ for some $A_j \subset G'$ and $g_j \in G$. Since $G'$ is a normal subgroup of $G$, we have, for each $i \in [1,t]$, that $A_i' = (A_i)^{(g_1 \ldots g_{i-1})^{-1}}$ for some $A'_i \subset G'$. It follows that $\pi(S_1)\ldots \pi(S_t) = (A_1g_1) \ldots (A_tg_t) = A_1'\ldots A_t'(g_1 \ldots g_t)$. In the special case where $G' \cong C_q$ with $q$ prime, then classical results on product-set cardinalities in $C_q$, such as the Cauchy-Davenport Theorem (Lemma \ref{lemCD}), may be applied to bound the cardinality of the product-set $\pi(S_1) \ldots \pi(S_t)$. Throughout the paper, this will be done without further reference to the intermediate sets $A_i'$.

\subsection{Ordered sequences over groups}

Let $\F^*(G)$ denote the free non-abelian monoid with basis $G$, that is, $\F^*(G)$ is the semigroup of words over the alphabet $G$. The elements of $\F^*$, called {\em ordered sequences} over $G$, are written as 
$$S^* = g_1 \bd \ldots \bd g_k = \prod_{j=1}^k g_j.$$

By disregarding the order of the elements in $\F^*(G)$, we obtain a natural map $[\cdot]: \F^*(G) \to \F(G)$. An ordered sequence $S^* \in \F^*(G)$ with $[S^*] = S$ is called an {\em ordering} of the sequence $S \in \F(G)$. Furthermore, if $S = [S^*]$, then we set $\supp(S^*) = \supp(S)$, $|S^*| = |S|$, and $\v_g(S^*) = \v_g(S)$ for every $g \in G$ to be the {\em support} of $S^*$, the {\em length} of $S^*$, and the {\em multiplicity} of $g \in G$ in $S^*$.

Let $S^* = g_1 \bd \ldots \bd g_k \in \F^*(G)$. For any subset $J \subset [1,k]$, set $S^*(J) = \prod_{j \in J} g_j$, where the order is taken in increasing order of the indices in $J$. We say that $S^*(J)$ is an {\em ordered subsequence} of $S^*$. For integers $0 \le i \le j$, we abbreviate $S^*(i,j) = S^*([i,j])$ and $S^*(j) = S^*(\{j\})$; the former is called a {\em consecutive subsequence}, while the latter denotes the {\em $j$-th term} of $S^*$. Moreover, $\pi: \F^*(G) \to G$ denotes the {\em product} of $S^*$ in the order the terms appear, that is, $\pi(g_1 \bd \ldots \bd g_k) = g_1 \ldots g_k$. If $S = [S^*]$, then it is clear that $\pi(S^*) \in \pi(S)$. A {\em factorization} of $S^* \in \F^*(G)$ (of length $t$) is a $t$-tuple $(S_1^*, \ldots, S_t^*)$ of nontrivial consecutive subsequences $S_i^* \mid S^*$ such that $S^* = S_1^* \bd \ldots \bd S_t^*$.

\subsection{On the group $C_q \rtimes C_p$}

We consider the groups $G$ of order $pq$, where $p \le q$ are prime numbers. If $p=q$, then either $G \cong C_{p^2}$ or $G \cong C_p^2$, whence $G$ is abelian. Suppose now that $p<q$. If $p \nmid q-1$, then an immediate consequence of Sylow's Theorem is that $G \cong C_{pq}$ is cyclic. It remains to analyse the case $p \mid q-1$. Another application of Sylow's Theorem yields, up to isomorphism, exactly two groups of order $pq$: the cyclic group $C_{pq}$, and only one non-abelian group, which can be written as the semidirect product $C_q \rtimes C_p$. The particular case $p=2$ corresponds to the dihedral group of order $2q$, which has been extensively studied (see \cite{Bas,MR1,GeGr,OhZh,OW}), therefore we will assume that $p,q$ are both odd prime numbers with $p \mid q-1$. From now on, we denote 
\begin{equation}\label{defG}
	\G = C_q \rtimes C_p \cong \langle \alpha, \tau \colon \alpha^q = \tau^p = 1, \alpha \tau = \tau \alpha^s \rangle,
\end{equation}
where $s$ has order $p$ modulo $q$. The commutator subgroup of $\G$ is $\G' = \langle \alpha \rangle$ and its center is $Z(\G) = \{1\}$. The centralizer of $g \in \G \backslash \{1\}$ is $C_{\G}(g) = \langle g \rangle$. Moreover, $\ord(g) = q$ for every $g \in \G' \backslash \{1\}$ and $\ord(g) = p$ for every $g \in \G \backslash \G'$. Since $p,q$ are odd and $p \mid q-1$, it follows that $q \ge 2p+1$.

It is worth mentioning that the direct and inverse problems over $\G$ associated to other invariants are already known, such as $\d(\G) = p+q-2$ \cite[Lemma 14]{Bas} (see also \cite{MR2} for the inverse problem), and the Erd\H os-Ginzburg-Ziv constant \cite[Theorem 15]{Bas} (see also \cite{QL1} for the inverse problem).

\section{Preliminary results}\label{seclemas}

The proof of Theorem \ref{main} closely follows the approach used by Grynkiewicz \cite{Gr} in his solution of the direct problem. Several lemmas from his paper and also from \cite{GeGr} are employed, some of them in an adapted form. We begin by stating a few general lemmas.

\begin{lem}[Cauchy-Davenport Theorem {\cite[Theorem 6.2]{Gr1}}]\label{lemCD}
	Let $G \cong C_q$, where $q$ is a prime number, and let $A, B \subset G$ be non-empty subsets. Then
	$$|AB| \ge \min\{q, |A|+|B|-1\}.$$
\end{lem}

\begin{lem}[{\cite[Lemma 2.1]{GeGr}}]\label{lem23} 
	Let $G$ be a group, let $U^*\in\mathcal{F}^*(G)$ be an ordered sequence with $\pi(U^*)=1$ and let  $[U^*]\in\mathcal{A}(G)$ an atom. Then there are no consecutive product-one subsequences of $U^*$ that are proper and nontrivial.
\end{lem}

\begin{lem}[{\cite[Lemma 2.2]{GeGr}}]\label{lem24} 
	Let $G$ be a group and let $S\in\mathcal{F}(G)$ be a product-one sequence. If $T\mid S$ is a subsequence with $\pi(T)\subseteq G'$, then $\pi(S\bd T^{[-1]})\subseteq G'$. In particular, if $T\mid S$ is a product-one subsequence, then $\pi(S\bd T^{[-1]})\subseteq G'$. 
\end{lem}

\begin{lem}[{\cite[Lemma 2.4.1]{GeGr}}]\label{lem26} 
	Let $G$ be a finite group. Then every ordered sequence $S\in\mathcal{F}^*(G)$ of length $|S|\ge|G|$ has a consecutive, product-one subsequence that is nontrivial. In particular, we have ${\sf{d}}(G)+1\le {\sf{D}}(G)\le|G|$.
\end{lem}

\begin{lem}[{\cite[Theorem 5.4.5.2]{GeHK}, see also \cite[Theorem 2.1]{GaGeSc}}]\label{lem5452}
	Let $G \cong C_n$, where $n \ge 3$, and let $S \in \F(G)$ be a product-one free sequence of length $|S| \ge \frac{n+1}{2}$. Then there exists $g \in \supp(S)$ such that $\v_g(S) \ge 2|S| - n + 1$. In particular, $\D(G) = n$ and if $|S| = n-1$, then $S = g^{[n-1]}$.
\end{lem}

The following lemma is a crucial technical tool that will be invoked repeatedly in the proof of the main theorem. It embodies a simple yet effective algorithm, and for further details on the underlying idea we refer to the discussion after Lemma 3.2 in \cite{Gr}.

\begin{lem}[{\cite[Lemma 3.3]{Gr}}]\label{lem33} 
	Let $G$ be a non-abelian finite group, let $S^*\in\mathcal{F}^*(G)$ be an ordered sequence, let $H\le G$ be an abelian subgroup, let $$\omega\ge1,\quad \omega_H\in\Z,\quad \mbox{ and }\quad \omega_0\in\{0\}\cup[2,|S^*|]\quad \mbox{ with }\quad \omega_0\le \omega,$$ and suppose that $|\pi(S_0)|\ge|S_0|=\omega_0$ and $\pi(S_0)\cap(G\backslash Z(G))\neq\varnothing$ (if $\omega_0>0$), where $S_0=[S^*(1,\omega_0)]$, and that there are at least $\omega_H$ terms of $S\bd S_0^{[-1]}$ from $H$.
	
	Then there exists an ordered sequence $S'^*\in\mathcal{F}(G)$ with 
	$$[S'^*]=[S^*]\quad \mbox{ and }\quad \pi(S'^*)\in\pi(S^*)^G$$ 
	having a factorization $$S'^*=T_1^*\bd\ldots\bd T_{r-1}^*\bd T_r^*\bd R^*,$$ where $T_1^*,\dots,T_r^*,R^*\in\mathcal{F}^*(G)$ and $r\ge0$, such that, letting $R=[R^*]$ and $T_i=[T_i^*]$ for $i\in[1,r]$, we have $S_0\mid T_1$ (if $\omega_0>0$), $$\pi(T_i)\cap(G\backslash Z(G))\neq\varnothing\quad \mbox{ and }\quad |\pi(T_i)|\ge|T_i|\ge2\quad \mbox{ for }i\in[1,r],\quad \pi(T_i)^G=\pi(T_i)\quad \mbox{ for }i\in[1,r-1],$$ and either \begin{enumerate}\item[(i)] $\sum_{i=1}^r|T_i|\le\omega-1$ and $\langle\supp(R)\rangle<G$ is a proper subgroup, or\item[(ii)] $\omega\le\sum_{i=1}^r|T_i|\le\omega+1$, with the upper bound only possible if $|T_r|=2$ and $\sum_{i=1}^{r-1}|T_i|=\omega-1$, and there are at least $\omega_H$ terms of $R$ from $H$, or\item[(iii)] $\sum_{i=1}^r|T_i|\le\omega-1$ and there are precisely $\omega_H$ terms of $R$ from $H$.\end{enumerate} 
\end{lem}

From now on, all results are restricted to the group $\G$ defined in Equation \eqref{defG}.


\begin{lem}[{\cite[Lemma 5.3]{Gr}}]\label{lem53} 
	Let $S\in\mathcal{F}(\G\backslash\{1\})$ and $g\in \G\backslash \G'$. Then $|\pi(g\bd S)|\ge\min\{q,|g\bd S|\}$. 
\end{lem}

\begin{lem}[{\cite[Lemma 5.4]{Gr}}]\label{lem54} 
	Let $S\in\mathcal{F}(\G'\backslash\{1\})$ and $g_1,g_2\in \G\backslash \G'$. Suposse $g_1g_2\notin \G'$. Then $$|\pi(g_1\bd g_2\bd S)|\ge\min\{q,2|S|+1\}.$$
\end{lem}

\begin{lem}[{\cite[Lemma 5.5]{Gr}}]\label{lem55} 
	Let $S\in\mathcal{F}(\G\backslash\{1\})$. If $\langle\supp(S)\rangle=\G$, then $|\pi(S)|\ge\min\{p,|S|\}$.
\end{lem}

\begin{lem}[{\cite[Lemma 5.8]{Gr}}]\label{lem58} 
	Let $S\in\mathcal{F}(\G)$. If $|S|\ge q+2p-3$, then there is a nontrivial, product-one subsequence $T\mid S$ with $|T|\le q$. 
\end{lem}

\begin{lem}[{\cite[Lemma 5.11]{Gr}}]\label{lem511} 
	Let $T_1,\dots,T_r\in\mathcal{F}(\G)$ be sequences for which $$\pi(T_i)\cap(\G\backslash Z(\G))\neq\varnothing\quad \mbox{ and }\quad |\pi(T_i)|\ge|T_i|\ge2\quad \mbox{ for }i\in[1,r],\quad \pi(T_i)^{\G}=\pi(T_i)\quad \mbox{ for }i\in[1,r-1].$$ 
	Then the following hold:
	\begin{enumerate}
		\item $|\pi(T_1)\ldots\pi(T_r)|\ge\min\left\{q-1,\sum_{i=1}^r|\pi(T_i)|\right\}\ge\min\left\{q-1,\sum_{i=1}^r|T_i|\right\};$
		\item if $\sum_{i=1}^r|T_i|\ge q+1$, then $|\pi(T_1)\ldots\pi(T_r)|=q.$
	\end{enumerate}
\end{lem}



The following result is essentially contained in the proof of \cite[Lemma 14]{Bas} (see also \cite[Lemma 4]{BRGS}), and we reproduce it here for convenience.

\begin{lem}[]\label{lem4} 
	If $T\in\mathcal{F}(\G)$ such that $|T|\ge q$ and $\pi(T)\cap \G'\neq\varnothing$, then $T$ contains a product-one subsequence.
\end{lem}

\proof 
	Since $\G/\G'\cong C_p$, which is an abelian group, we obtain that $\pi(T)\subseteq \G'$. Let us factorize $T$ as 
	$$T = T_1 \bd \dots \bd T_k,$$ 
	where $\pi(T_i)\subseteq \G'$ and each $T_i$ is minimal with respect to this property. Since $|T| \ge q > p = \d(\G/\G') + 1$ and $2p + 1 \le q$, we obtain that $k \ge 3$. Let us denote $A_1 = \pi(T_1)$ and $A_i = \pi(T_i)\cup\{1\}$ for all $i \in [2, k]$. 

	We claim that either $1 \in \pi(T_i)$ for some $i \in [1, k]$, in which case we conclude the proof, or $1 \notin \pi(T_i)$ and $|\pi(T_i)| \ge |T_i|$ for all $i \in [1, k]$.  In the latter case, it follows from the Cauchy–Davenport Theorem (Lemma \ref{lemCD}) that $$|A_1 \ldots A_k| \ge \min\left\{q, \sum_{i=1}^k |A_i| - (k - 1)\right\}= \min\left\{q, \sum_{i=1}^k |\pi(T_i)| + (k - 1) - (k - 1)\right\}\ge \min\left\{q, \sum_{i=1}^k |T_i|\right\} = q. $$ Since $A_1 \ldots A_k \subseteq \Pi(T)$, we conclude that $T$ has a product-one subsequence.

	We remark that if $g_1,g_2 \in \G$ are such that $g_1g_2 = g_2g_1 \in \G'\backslash\{1\}$, then $g_1,g_2 \in \G'\backslash\{1\}$. Let $g_1 = \uptau^{a_1}\alpha^{b_1}$ and $g_2 = \uptau^{a_2}\alpha^{b_2}$. Hence, $$\uptau^{a_1 +a_2} \alpha^{b_1 + b_2 s^{a_1}} = \uptau^{a_1 + a_2} \alpha^{b_2 + b_1 s^{a_2}} = \alpha^u \in \G'\backslash\{1\},$$ and we have $a_1 +a_2 \equiv 0 \pmod{p}$ and $b_1 + b_2 s^{a_1} \equiv b_2 + b_1 s^{a_2} \equiv u \pmod{q}$. Then $$u \equiv b_1 + b_2 s^{a_1} \equiv s^{a_1}(b_1 s^{-a_1} + b_2) \equiv s^{a_1}(b_1 s^{a_2} + b_2) \equiv u s^{a_1} \pmod{q}.$$ Since $u \not\equiv 0 \pmod{q}$, we obtain that $a_1 \equiv 0 \equiv a_2 \pmod{p}$, and thus $g_1,g_2 \in \G'\backslash\{1\}$.

	Now we finish the proof by showing that if $1 \notin \pi(T_i)$, then $|\pi(T_i)| \ge |T_i|$. Let $|T_i| = t_i$, write $T_i = g_1 \bd \dots \bd g_{t_i}$, and define $\pi_j(T_i) = g_j \ldots g_{t_i} g_1 \ldots g_{j-1}$ for $j \in [1, t_i]$. Observe that $|\{\pi_j(T_i)\}| \le |\pi(T_i)|$. We will show that $|\{\pi_j(T_i)\}| = t_i$, that is, $|T_i| \le |\pi(T_i)|$. If $|\{\pi_j(T_i)\}| < t_i$, then there exist $1 \le k < \ell \le t_i$ such that $\pi_k(T_i) = \pi_\ell(T_i)$. Then 
	$$g_k \ldots g_{t_i} g_1 \ldots g_{k-1} = g_\ell \ldots g_{t_i} g_1 \ldots g_{\ell-1}.$$ 
	Let $g = g_k \ldots g_{\ell-1}$ and $g' = g_\ell \ldots g_{t_i} g_1 \ldots g_{k-1}$. Hence, $gg' = g'g \in \G'\backslash\{1\}$, since $1 \notin \pi(T_i)$. But it follows from the remark above that $g,g' \in \G'\backslash\{1\}$, which contradicts the minimality of $T_i$. Therefore, $|\{\pi_j(T_i)\}| = t_i$, and we conclude the proof.

\qed

The following results are adapted from \cite{Gr}. Their proofs are similar.

\begin{lem}[{Adapted from \cite[Lemma 5.9]{Gr}}]\label{lem59} 
	Let $S\in\A(\G)$ with $\v_{\G\backslash \G'}(S)\ge3$. If $|S|\ge2q$, then 
	$$\v_{\G'}(S) \le \dfrac{q-3}{2}.$$
\end{lem}

\proof
	At first, we claim that there exist $g_1,g_2\in\supp(S)\cap(\G\backslash \G')$ such that $g_1g_2\notin \G'$. In fact, since $\v_{\G\backslash \G'}(S)\ge3$, let $x,y,z\in\supp(S)\cap(\G\backslash \G')$ and assume that $xy,xz,yx\in \G'$. Let $\phi_{\G'}(g)=g\G'$ be the canonical homomorphism. Then $\G'=\phi_{\G'}(xy)=\phi_{\G'}(xz)$ and, hence, $\phi_{\G'}(y)=\phi_{\G'}(z)$. This implies that $\G'=\phi_{\G'}(yz)=\phi_{\G'}(y)^2$. Since $|\G/\G'|=p$ is an odd prime number, we obtain that $y\in \G'$, which is a contradiction.

	Now, we will show that $\v_{\G'}(S)\le\frac{q-3}{2}$ and, to this end, let us assume otherwise, that is, assume that $\v_{\G'}(S)\ge\frac{q-1}{2}$. Let $T\mid S$ be a subsequence such that $\supp(T)\subseteq \G'$ and $|T|=\frac{q-1}{2}$. Then 
	$$|S\bd(g_1\bd g_2\bd T)^{[-1]}|=|S|-|T|-2\ge2q-\dfrac{q-1}{2}-2=q+\dfrac{q-1}{2}-1\ge q+p-1={\sf{d}}(\G)+1.$$ 
	Let $R\mid S$ be a nontrivial product-one subsequence such that $g_1\bd g_2\bd T\mid S\bd R^{[-1]}$. It follows from Lemma \ref{lem24} that $\pi(S\bd R^{[-1]})\subseteq \G'$. However, since $\supp(T)\subseteq\supp(S)\subseteq \G\backslash\{1\}$,  it follows from Lemma \ref{lem54} that $$|\pi(g_1\bd g_2\bd T)|\ge\min\left\{q,2|T|+1\right\}=q,$$ which implies that $1\in\pi(S\bd R^{[-1]})$. Therefore, $R\bd(S\bd R^{[-1]})$ is a nontrivial factorization of $S$ into two product-one subsequences, contradicting the fact that $S$ is an atom.
\qed

\begin{lem}[{Adapted from \cite[Lemma 5.12]{Gr}}]\label{lem512} 
	Let $S\in\mathcal{A}(\G)$. If $|S|\ge2q$, then 
	$$\v_H(S)\le q-1\quad \mbox{ for every subgroup }\quad H\le \G\quad \mbox{ with }\quad |H|=p.$$
\end{lem}

\proof
	Let $S^* \in \mathcal{F}(\G)$ be such that $\pi(S^*) = 1$ and $S = [S^*]$, and let us assume that there exists a subgroup $H \leq \G$ with $|H| = p$ such that $\v_H(S) \geq q$. We will apply Lemma \ref{lem33} to $S^*$ using $H$ with  
	$$\omega = q + 1, \quad \omega_H = p + 1, \quad \mbox{and} \quad \omega_0 = 0.$$  
	Let  
	$$S'^* = T_1^* \bd \ldots \bd T_r^* \bd R^*$$  
	be the factorization obtained in Lemma \ref{lem33}, and we will analyze the three cases. First, we remark that since $\pi(S'^*) \in \pi(S^*)^{\G} = \{1\}$, it follows that $\pi(S'^*) = 1$.  

	\vspace{2mm}
	\noindent
	{\bf Case 1.} $\sum_{i=1}^r |T_i| \leq \omega - 1 = q$ and $K = \langle \supp(R) \rangle < \G$ is a proper subgroup.  

		Since $\v_H(S) \geq q$ and $\sum_{i=1}^r |T_i| \leq q$, we have two possibilities: either there is at least one element of $H\backslash\{1\}$ in $K$, or there is none different from 1. 

		Assume first that there is a element of $H$ in $K$ and, since $K$ is a proper subgroup, we must have $H = K$.  However, since $|R| = |S| - \sum_{i=1}^r |T_i| \geq 2q - q = q > p$, Lemma \ref{lem26} guarantees a nontrivial and proper product-one consecutive subsequence of $R^*$, contradicting Lemma \ref{lem23} and the fact that $S$ is an atom.  

		Now, assume that $K$ contains no element of $H$ other than $1$. In this case, we have $\v_H(S) = q$, $\sum_{i=1}^r |T_i| = q$, and $|R| \ge q$. Let us denote $K = \langle \uptau^a \alpha^b \rangle$. If $a \not\equiv 0 \pmod{p}$, then $|K| = p < q \le |R|$. On the other hand, if $a \equiv 0 \pmod{p}$, that is, $K = \G'$, then ${\sf d}(\G') = q - 1 < q \le |R|$. In both cases, we also obtain a contradiction by Lemma \ref{lem26}, Lemma \ref{lem23}, and the fact that $S$ is an atom.

	\vspace{2mm}
	\noindent
	{\bf Case 2.}  $q+1 = \omega \le \sum_{i=1}^r |T_i| \le \omega+1 = q+2$ and there are at least $\omega_H=p+1$ elements of $R$ from $H$. 

		Since ${\sf{d}}(H)=p-1$, there exists a nontrivial product-one subsequence $R' \mid R$. As a consequence of Lemma \ref{lem24}, we have that $\pi(S \bd R'^{[-1]}) \subseteq \G'$. Now, observe that $T_1 \bd \ldots \bd T_r \mid S \bd R'^{[-1]}$ and, by Lemma \ref{lem511}.(2), $|\pi(T_1) \ldots \pi(T_r)| = q$. Then $\pi(S \bd R'^{[-1]}) = \G'$. Therefore, $S = R' \bd (S \bd R'^{[-1]})$ is a nontrivial factorization of $S$ into two product-one subsequences, contradicting the fact that $S$ is an atom.

	\vspace{2mm}
	\noindent
	{\bf Case 3.} $\sum_{i=1}^r|T_i|\le\omega-1=q$ and $\v_H(R)=\omega_H=p+1$.

		Since $\v_H(R) = \omega_H = p+1 < q \leq \v_H(S)$, we obtain that 
		$$\v_H(T_1 \bd \ldots \bd T_r) = \v_H(S) - \v_H(R) \geq q - p - 1.$$ 
		Since $H$ is abelian and $|\pi(T_i)| \geq |T_i| \geq 2$, each $T_i$ has terms from $\G \backslash H$. Then 
		$$\sum_{i=1}^r |T_i| \geq q - p - 1 + r \ge q - p,$$ 
		where $r \geq 1$. Notice that if $\sum_{i=1}^r|T_i|=q-p$, then $r=1$ and $|R|=q+p$. Now, since $\sum_{i=1}^r|T_i|\le q$, we obtain that
		\begin{equation*}
			\v_{\G\backslash H}(R)=
			\begin{cases}
				|S|-|T_1|-\v_H(R) = 2q-q+p-p-1=q-1, & \text{if } r=1,\\
				|S|-\sum_{i=1}^r|T_i|-\v_H(R)\geq 2q-q-p-1=q-p-1>p-1, & \text{if } r\ge 2.
			\end{cases}
		\end{equation*}
		Then there exists $R'\mid R$ such that
		\begin{equation*}
			|R'|=
			\begin{cases}
				q-1, & \text{if } r=1,\\
				p-1, & \text{if } r\ge2,
			\end{cases}
		\end{equation*}
		and $\operatorname{supp}(R')\cap H=\varnothing$. Let $g \in \operatorname{supp}(R) \cap H$. Then $\langle \supp(g \bd R') \rangle = \G$, and Lemma \ref{lem55} guarantees that 
		$$|\pi(g \bd R')| \geq \min\left\{p, |g \bd R'|\right\} = p.$$ 
		Since $\pi(g\bd R')$ is contained in a $\G'$-coset, we can apply the Cauchy-Davenport Theorem (Lemma \ref{lemCD}) and, together with Lemma \ref{lem511}.(1), we obtain that
		\begin{eqnarray*} 
			|\pi(T_1 \bd \ldots \bd T_r) \pi(g \bd R')| 
			&\geq& \min\left\{q, |\pi(T_1 \bd \ldots \bd T_r)| + |\pi(g \bd R')| - 1\right\}\\
			&\geq& \min\left\{q, \min\{q-1, \sum_{i=1}^r |T_i|\} + |\pi(g \bd R')| - 1\right\}
		\end{eqnarray*}
		and, since $q\ge2p+1$, 
		$$\min\left\{q, \min\left\{q-1, \sum_{i=1}^r |T_i|\right\} + |\pi(g \bd R')| - 1\right\}\ge
		\begin{cases}
    		\min\left\{q, \min\left\{q-1, q - p\right\} + q - 1\right\}=q&\mbox{ if }r=1,\\
    		\min\left\{q, \min\left\{q-1, q - p+1\right\} + p - 1\right\} = q&\mbox{ if }r\ge2.
		\end{cases}$$

		In both cases, since $\v_H(R) = p+1$ and $\v_H(g \bd R') = 1$, we still have $p = |H|$ terms of $R \bd (g \bd R')^{[-1]}$ from $H$. Since ${\sf{d}}(H) = p-1$, there exists a nontrivial product-one subsequence $R'' \mid R \bd (g \bd R')^{[-1]}$. It follows from Lemma \ref{lem24} that $\pi(S \bd R''^{[-1]}) \subseteq \G'$. However, since $T_1 \bd \ldots \bd T_r \bd g \bd R' \mid S \bd R''^{[-1]}$, it follows that $\pi(S \bd R''^{[-1]}) = \G'$. Therefore, $S = R'' \bd (S \bd R''^{[-1]})$ is a nontrivial factorization of $S$ into two product-one subsequences.
	
	\vspace{2mm}

	Summing up, we have that $\v_H(S)\le q-1$ for every subgroup $H\le \G$ with $|H|=p.$
	
\qed

\section{Proof of Theorem \ref{main}}\label{secprovamain}

The proof of Theorem \ref{main} follows from the next two results. The first states that any atom of length $\D(\G)$ must have at most two terms from $\G \backslash \G'$. 

\begin{teo}\label{casov>2}
	Let $S\in\mathcal{F}(\G)$ such that $|S|=2q$. If $\v_{\G\backslash \G'}(S)\ge3$, then $S\notin\mathcal{A}(\G)$.
\end{teo}

\proof 
	Let us assume that $S\in\mathcal{A}(\G)$ and let $S^*\in\mathcal{F}(\G)$ be such that $\pi(S^*)=1$ and $S=[S^*]$. It follows from Lemma \ref{lem59} that $\v_{\G'}(S)\le\frac{q-3}{2}$. In light of Lemma \ref{lem58}, since $|S|=2q>q+2p-3$, there exists $U\mid S$ such that $|U|\le q$ and $1\in\pi(U)$. Let $U$ be minimal with respect to this property, that is, $U$ is the shortest product-one subsequence of $S$. We will split the proof into three cases.

	\vspace{2mm}
	\noindent
	{\bf Case 1.} $|U|=q$.
		
		Let $V = S \bd U^{[-1]}$. It follows from Lemma \ref{lem24} that $\pi(V) \subseteq \G'$. Let us assume that $1 \notin \pi(V)$; otherwise, $S = U \bd V$ is a factorization of $S$ into two product-one subsequences, contradicting the fact that $S \in \mathcal{A}(\G)$. As a consequence of Lemma \ref{lem4}, $V$ contains a product-one subsequence of lenght at most $q - 1$, contradicting the minimality of $U$. 

	\vspace{2mm}
	\noindent
	{\bf Case 2.} $|U|\le q-p$.

		We first claim that $U$ can be taken as a nontrivial product-one subsequence with $|U|\le p$ and 
		$|\langle\supp(U)\rangle|=p$. Let $W=S\bd U^{[-1]}$. If $\v_{\G'}(W)=0$, let $W_0$ be the trivial sequence. Otherwise, since 
		$$|W|=|S|-|U|\ge2q-q+p=q+p>\dfrac{q-3}{2}\ge\v_{\G'}(S),$$ let $W_0\mid W$ be the subsequence containing all terms of $W$ from $\G'$ and exactly one term from $\G\backslash \G'$. Moreover, if $W_0$ is nontrivial, Lemma \ref{lem53} guarantees that $|\pi(W_0)|\ge|W_0|\ge2$ and, hence, $\pi(W_0)\cap(\G\backslash\{1\})\neq\varnothing$. Let us fix $W^*$ as any ordering of $W$ such that $[W^*(1,|W_0|)]=W_0$. Now, we apply Lemma \ref{lem33} to $W$ by taking $H=\{1\}$, $\omega=q+1$, $\omega_H=-1$, and $\omega_0=|W_0|\le\frac{q-1}{2}$. Let $W'^*=T_1^*\bd\ldots\bd T_r^*\bd R^*$ be the factorization obtained in Lemma \ref{lem33}. We note that since $\omega_H=-1$, the third case of this lemma does not occur. If the second case holds, that is, $q+1\le\sum_{i=1}^r|T_i|$, then Lemma \ref{lem511}.(2) implies that $|\pi(T_1\bd\ldots\bd T_r)|\ge|\pi(T_1)\ldots\pi(T_r)|=q$. However, since $1\in\pi(U)$, then $\pi(S\bd U^{[-1]})=\pi(W)\subseteq \G'$ by Lemma \ref{lem24}. Thus, $\pi(W)=\G'$, and hence, $S=W\bd U$ is a nontrivial factorization of $S$ into two product-one subsequences, contradicting the fact that $S$ is an atom. Finally, assume that the first case of Lemma \ref{lem33} holds, that is, $\sum_{i=1}^r|T_i|\le\omega-1=q$ and $K=\langle\supp(R)\rangle$ is a proper subgroup of $\G$. In this case, $$|R|=|W|-\sum_{i=1}^r|T_i|\ge q+p-q=p.$$ Since $W_0\mid T_1$, there is no term in $R$ from $\G'$ and thus $|K|=p$. Then $|R|\ge p={\sf{d}}(K)+1$, and there exists a nontrivial product-one subsequence of $R$ with at most $p$ elements. This proves that there exists a nontrivial product-one subsequence $U$ of $S$ with $|U|\le p$ and $|\langle\supp(U)\rangle|=p$. From now on, let $U$ be this product-one subsequence of $S$. 

		Let us define $W=S\bd U^{[-1]}$, $W_0$ and $W^*$ as done above. We will use Lemma \ref{lem33} on $W$ again, using the same parameters: $H=\{1\}$, $\omega=q+1$, $\omega_H=-1$, and $\omega_0=|W_0|\le\frac{q-1}{2}$. Statements (2) and (3) of Lemma \ref{lem33} do not hold, as argued in the paragraph above. Then let us assume that $\sum_{i=1}^r|T_i|\le \omega-1=q$ and that $K=\langle\supp(R)\rangle$ is a proper subgroup of $\G$. Then
		$$|R|=|S|-|U|-\sum_{i=1}^r|T_i|\ge 2q-p-q=q-p\ge p+1.$$ 
		Let $K'=\langle\supp(U)\rangle$. If $K=K'$, then all terms of $R\bd U$ belong to the same subgroup of order $p$, and $|R\bd U|=|S|-\sum_{i=1}^r|T_i|\ge 2q-q=q$, which contradicts Lemma \ref{lem512}. Therefore, $K\neq K'$. Since $|R|>{\sf{d}}(K)$, there exists a nontrivial product-one subsequence $L\mid R$ with $\langle\supp(L)\rangle=K$.

		Let us define $V=S\bd L^{[-1]}=W\bd U\bd L^{[-1]}$ and $Z=R\bd U\bd L^{[-1]}\mid V$. First, we observe that since $1\in\pi(L)$, then $\pi(V)\subseteq \G'$ by Lemma \ref{lem24}. Moreover, there are terms of $Z$ from both $K$ and $K'$. Since $1\notin\supp(S)$, there exist elements $g\in K\backslash\{1\}$ and $g'\in K'\backslash\{1\}$ such that $g,g'\in\supp(Z)$. Since $K\neq K'$, we obtain that $gg'\neq g'g$. If $q-1\le\sum_{i=1}^r|T_i|\le q$, since $\pi(g\bd g')$ is contained in a $\G'$-coset, we apply the Cauchy-Davenport Theorem (Lemma \ref{lemCD}) and, together with Lemma \ref{lem511}.(1), we obtain that
		\begin{eqnarray*} 
 			|\pi(T_1\bd\ldots\bd T_r)\pi(g \bd g')| 
			&\geq& \min\left\{q, |\pi(T_1\bd\ldots\bd T_r)| + |\pi(g \bd g')| - 1\right\}\\
			&\geq& \min\left\{q, \min\left\{q-1, \sum_{i=1}^r |T_i|\right\} + |\pi(g \bd g')| -1\right\}\\
			&\geq& \min\left\{q, \min\left\{q-1, q - 1\right\} + 2 -1\right\} = q.
		\end{eqnarray*}
		Since $T_1\bd\ldots\bd T_r\bd g\bd g'\mid V$, we obtain that $\pi(V)=\G'$, and hence, $S=V\bd L$ is a nontrivial factorization of $S$ into two product-one subsequences. Therefore, $\sum_{i=1}^r|T_i|\le q-2$.

		Let us fix $V_0=T_1\bd\ldots\bd T_r$ and $V^*$ as any ordering of $V$ such that $[V^*(1,|V_0|)]=V_0$. Then $|V_0|\le q-2$, and as a consequence of Lemma \ref{lem511}.(1),
		$$|\pi(V_0)|=|\pi(T_1\bd\ldots\bd T_r)|\ge|\pi(T_1)\ldots\pi(T_r)|\ge\min\left\{q-1,\sum_{i=1}^r|T_i|\right\}=|V_0|.$$ 
		Let us apply Lemma \ref{lem33} again but with the parameters $H=\{1\}$, $\omega=q+1$, $\omega_H=-1$, and $\omega_0=|V_0|\le q-2$. Again, the statements (2) and (3) from Lemma \ref{lem33} do not hold. Let $V'^*=T_1'^*\bd\ldots\bd T_r'^*\bd R'$ be the factorization obtained from Lemma \ref{lem33}. Since $V_0\mid T_1'$, then $T_1\bd\ldots\bd T_r\mid T_1'$ and $R'\mid Z$. Now, $\supp(Z)\subseteq K\cup K'$ with $\v_{K'}(Z)=\v_{K'}(U)=|U|\leq p$. Then at most $p$ terms of $R'$ are from $K'$, and the remaining are from $K$. However, $\langle\supp(R')\rangle$ is a proper subgroup of $\G$ and
		$$|R'|=|V'|-\sum_{i=1}^r|T_i'|=|S|-|L|-\sum_{i=1}^r|T_i'|\ge 2q-p-q=q-p\ge p+1.$$ 
		Thus, all terms of $R'$ are from $K$. But, since $\supp(R'\bd L)\subseteq K$ and
		$$|R'\bd L|=|S|-\sum_{i=1}^r|T_i'|\ge 2q-q=q,$$
		we obtain a contradiction by Lemma \ref{lem512}. Therefore, this case cannot occur.

	\vspace{2mm}
	\noindent
	{\bf Case 3.} $q-p<|U|\le q-1$.

		Since $|S|-|U|\ge q+1>\frac{q-3}{2}+p$, there exist $g_1,g_2\in\supp(S\bd U^{[-1]})\cap(\G\backslash \G')$ with $\langle g_1,g_2\rangle=\G$. Let us fix $W=S\bd(U\bd g_1\bd g_2)^{[-1]}$. Since $|W|=|S|-|U|-2\ge q-1>\v_{\G'}(S)$, let us proceed as in case 2. and let $W_0$ be the trivial sequence if $\v_{\G'}(W)=0$, and let $W_0$ be the sequence consisting of all terms of $W$ from $\G'$ and one term from $\G\backslash \G'$. In the latter case, Lemma \ref{lem53} guarantees that $|\pi(W_0)|\ge|W_0|\ge2$ and, hence, $\pi(W_0)\cap(\G\backslash\{1\})\neq\varnothing$. Let $W^*\in\mathcal{F}(\G)$ be any ordering of $W$ such that $[W^*(1,|W_0|)]=W_0$, and let us apply Lemma \ref{lem33} with the parameters $H=\{1\}$, $\omega={q-p-1}$, $\omega_H=-1$, and $\omega_0=|W_0|\le\frac{q-1}{2}\le\omega$. Let $W'^*=T_1^*\bd\ldots\bd T_r^*\bd R^*$ be the factorization obtained in Lemma \ref{lem33}. Since $\omega_H=-1$, statement (3) of this lemma does not occur. Let us analyze the other two statements.
		
		\vspace{2mm}
		\noindent
		{\bf Subcase 3.1.} Assume that the second statement of Lemma \ref{lem33} holds, that is, $$\omega=q-p-1\le\sum_{i=1}^r|T_i|\le\omega+1=q-p<q-1\le|W|.$$  
			As a consequence of Lemma \ref{lem511}.(1): if $\sum_{i=1}^r|T_i|=q-p-1$, then $$|R|=|W|-\sum_{i=1}^r|T_i|\ge (q-1)-(q-p-1)=p$$ and $$|\pi(T_1)\ldots\pi(T_r)|\ge \min\left\{q-1,\sum_{i=1}^r|T_i|\right\}\ge\min\left\{q-1,q-p-1\right\}=q-p-1;$$
			if $\sum_{i=1}^r|T_i|=q-p,$ then $$|R|=|W|-\sum_{i=1}^r|T_i|\ge (q-1)-(q-p)=p-1$$ and $$|\pi(T_1)\ldots\pi(T_r)|\ge \min\left\{q-1,\sum_{i=1}^r|T_i|\right\}\ge\min\left\{q-1,q-p\right\}=q-p.$$
			That is, $$|R|\ge p-1+\epsilon$$ and $$|\pi(T_1)\ldots\pi(T_r)|\ge q-p-\epsilon,$$ where $\epsilon\in\{0,1\}$. 
			As a consequence of Lemma \ref{lem53}, $$|\pi(R\bd g_1\bd g_2)|\ge\min\left\{q,|R|+2\right\}=p+1+\epsilon.$$ Now, using the Cauchy-Davenport Theorem (Lemma \ref{lemCD}), we obtain that $$|(\pi(T_1)\ldots\pi(T_t))(\pi(R\bd g_1\bd g_2))|\ge\min\left\{q, |\pi(T_1)\ldots\pi(T_t)|+|\pi(R\bd g_1\bd g_2)|-1\right\}=q.$$
			Then $1\in\pi(S\bd U^{[-1]})=\pi(W\bd g_1\bd g_2)$ and $S\notin\mathcal{A}(\G)$.

		\vspace{2mm}
		\noindent
		{\bf Subcase 3.2.} Assume now that the first statement of Lemma \ref{lem33} holds, that is, $\sum_{i=1}^r|T_i|\le\omega-1=q-p-2$ and $\langle{\supp}(R)\rangle<\G$ is a proper subgroup. Then $$|R|=|W|-\sum_{i=1}^r|T_i|\ge (q-1)-(q-p-2)=p+1.$$ Since $\v_{\G'}(R)=0$, $\langle{\supp}(R)\rangle$ must have order $p$. Since $|R|>{\sf d}(\langle{\supp}(R)\rangle)+1=p$, there exists $L\mid R$ such that $1\in\pi(L)$ and $|L|\le p$. But $|L|\le p<q-p<|U|$, which contradicts the minimality of the length of $U$, and this case cannot hold. This completes the proof.
		
\qed

In order to describe the minimal product-one sequences of maximal length $\D(\G)$, we need to study product-one sequences $S$ for which $\v_{\G\backslash \G'}(S) = 2$.

\begin{prop}\label{inv}
	Let $S \in \mathcal A(\G)$ with $|S| = 2q$ and $|S_{\G \backslash \G'}| = 2$. Then there exist $x,y \in \G$ and $s \in \Z_q^*$ for which $\G = \langle x,y \mid x^p = y^q = 1, yx = xy^s, \ord_q(s)=p \rangle$ and 
    \begin{equation}\label{formadeS}
		S = y^{[q-1]} \bd x \bd y^{[q-1]} \bd x^{p-1}y^{s^{p-1}+1}.
	\end{equation}
\end{prop}

\proof
	It is easy to verify that the sequence $S$ given by Eq. \eqref{formadeS} is a product-one sequence, since 
	$$y^{q-1} \cdot x \cdot y^{q-1} \cdot x^{p-1}y^{s^{p-1}+1} = y^{-1} \cdot x \cdot x^{p-1} \cdot y^{-s^{p-1}} \cdot y^{s^{p-1}+1} = 1.$$ 
	We claim that $S$ is minimal. Indeed, if $S$ is not minimal, then $S = S_1 \bd S_2$ with $S_1, S_2 \in \mathcal F(\G)$ both nontrivial product-one sequences. It is clear that either $S_1 = y^{[q]}$ or $S_2 = y^{[q]}$. Say $S_1 = y^{[q]}$. This implies that $S_2 = y^{[q-2]} \bd x \bd x^{p-1}y^{s^{p-1}+1}$. Since $S_2$ has product-one, it follows that 
	$$1 = y^{t} \cdot x \cdot y^{q-2-t} \cdot x^{p-1}y^{s^{p-1}+1} = y^t \cdot x \cdot x^{p-1} \cdot y^{(-2-t)s^{p-1}+s^{p-1}+1} = y^{t+(-2-t)s^{p-1}+s^{p-1}+1}$$
	for some $t \in [0,q-2]$. Therefore $s^{p-1}(t+1) - (t+1) \equiv 0 \pmod q$. Since $1 \le t+1 \le q-1$, we have $\gcd(t+1,q) = 1$ and hence $s^{p-1} \equiv 1 \pmod q$, a contradiction since $\ord_q(s) = p$.

	On the other hand, let $S$ be a minimal product-one sequence of length $2q$ for which $|S_{\G \backslash \G'}| = 2$. We may write 
	$$S = x^{a}y^{b_1} \bd x^{p-a}y^{b_2} \bd \prod_{i=1}^{2q-2} y^{c_i}, \quad a \in [1,p-1], \;\; b_1, b_2 \in [0,q-1], \;\; c_i \in [1,q-1].$$

	We may assume that $a = 1$ and $b_1 = 0$. Indeed,  
	$$\begin{cases}
		\ord_q(s^a) = p, \\
		(x^{a}y^{b_1})^p = x^{a p}y^{b_1 + b_1 s^{a} + b_1 s^{2a} + \ldots + b_1 s^{(p-1)a}} = y^{b_1 \left( \frac{s^{a p} - 1}{s^{a}-1} \right)} = 1, \\
		y \cdot x^a y^{b_1} = x^a y^{s^a + b_1} = x^a y^{b_1} \cdot y^{s^a},
	\end{cases}$$
	and this means that $\{x^a y^{b_1},y\}$ and $\{x, y\}$ generate non-abelian groups of order $pq$, which is unique up to isomorphism. Therefore both $\{x^a y^{b_1},y\}$ and $\{x, y\}$ generate isomorphic groups. 

	Since $S$ is a product-one sequence, it follows that $$\pi^*(T_1) \cdot x \cdot \pi^*(T_2) \cdot x^{p-1}y^{b_2} = 1, \quad T_1, T_2 \in \mathcal F(\G').$$
	If either $T_1$ or $T_2$ is a product-one sequence, then $S$ is not an atom. By Lemma \ref{lem5452}, we must have $T_j = (y^{c_j})^{[q-1]}$ for some $c_j \in [1,q-1]$, $j = 1,2$. Similarly to the previous paragraph, we may assume that $c_1 = 1$ since $\{x,y^{c_1}\}$ and $\{x,y\}$ generate isomorphic groups. Therefore 
	$$S = x \bd x^{p-1}y^{b_2} \bd y^{[q-1]} \bd (y^{c_2})^{[q-1]},$$ 
	so that 
	$$x^{p-1}y^{b_2}\cdot y^{q-1}\cdot x \cdot y^{c_2(q-1)} = y^{(b_2-1)s-c_2} = 1$$ 
	if and only if 
	$$(b_2-1)s \equiv c_2 \pmod q.$$
	In this case, if $c_2 \neq 1$, then, by Lemma \ref{lem5452}, $y^{[c_2]}\bd(y^{c_2})^{[q-2]} \mid T_1\bd T_2$ has a nontrivial ordered product-one subsequence $T_0^*=y^{[\ell]}\bd(y^{c_2})^{[k]}$, with $1\le\ell\le c_2$, $1\le k\le q-2$ and $\ell+c_2k\equiv0\pmod q$, which can be chosen to be consecutive, so that 
	$$1 = y^{(b_2-1)s-c_2} =  x^{p-1}y^{b_2}\cdot y^{c_2}\cdot y^{q-1-c_2}\cdot x\cdot (y^{c_2})^{q-2-k}\cdot y^{c_2-\ell} \in \pi(S \bd T_0^{[-1]}).$$ 
	Therefore $T_0 \bd (S \bd  T_0^{[-1]})$ is a decomposition of $S$ into nontrivial product-one subsequences. Thus $c_2 = 1$ and this completes the proof. 
	 
\qed

\vspace{1mm}

Now we are able to prove the main theorem of this paper.

\vspace{1mm}

\proofthmmain 
	Let $S \in \A(\G)$ with $|S| = 2q$. By Theorem \ref{casov>2}, $\v_{\G\backslash \G'}(S) \le 2$. If $\v_{\G\backslash \G'}(S) = 0$, then $S \in \F(\G')$. Since $\G' \cong C_q$, it follows that $\D(\G') = q$, therefore $S$ is not an atom. If $\v_{\G\backslash \G'}(S) = 1$, then $\pi(S) \cap \G' = \varnothing$, therefore $S$ is not a product-one sequence. This implies that $\v_{\G\backslash \G'}(S) = 2$. By Proposition \ref{inv}, it follows that $S$ is of the form \eqref{formaS}, and we are done.
\qed

\section{The union of sets of lengths containing $k$ and the $k$-th elasticity of $\B(C_q \rtimes C_p)$}\label{secaplica}

In this section, a {\em monoid} is a commutative cancelative semigroup with unit element. Suppose that $\HH$ is an atomic monoid, that is, every non-unit element can be written as a finite product of atoms, and let $\A(\HH)$ denote the set of atoms (irreducible elements) of $\HH$. In this sense, if $a \in \HH$, then there exist $u_1, \dots, u_k \in \A(\HH)$ such that $a = u_1 \dots u_k$. This $k$ is called the {\em length of the factorization} of $a$, and we define the {\em set of lengths} of $a$ as 
$$\L(a) = \{k \in \N \colon a \text{ has a factorization of length } k\}.$$ 
The {\em system of sets of lengths} of $\HH$ is $$\LL(\HH) = \{\L(a) \colon a \in \HH\}.$$ 
If not every element of $\HH$ is invertible, then, for $k \in \N$, the {\em union of sets of lengths containing $k$} is 
$$\U_k(\HH) = \bigcup_{L \in \LL(\HH) \atop k \in L} L.$$
Let $\rho_k(\HH) = \sup \U_k(\HH)$ be {\em $k$-th elasticity of $\HH$}, and let $\lambda_k(\HH) = \inf \U_k(\HH)$. For a subset $L \subset \N$, let $\rho(L) = \frac{\sup L}{\min L} \in \Q_{\ge 1} \cup \{\infty\}$ be the {\em elasticity of $L$}. The {\em elasticity of $\HH$} is defined as $\rho(\HH) = \sup\{\rho(L) \colon L \in \LL(\HH)\}$. It is possible to show that 
$$\rho(\HH) = \sup\left\{ \frac{\rho_k(\HH)}{k} \colon k \in \N \right\} = \lim_k \frac{\rho_k(\HH)}{k} \quad \text{ and } \quad  \frac{1}{\rho(\HH)} = \inf\left\{ \frac{\lambda_k(\HH)}{k} \colon k \in \N \right\} = \lim_k \frac{\lambda_k(\HH)}{k}$$(see \cite[Proposition 2.4]{Ge}). We have that 
$$\U_k(\HH) = \{\ell \in \N \colon \text{ there exist $u_1, \dots, u_k, v_1, \dots, v_\ell \in \A(\HH)$ such that $u_1 \dots u_k = v_1 \dots v_\ell$}\}.$$
From this, it is clear that $k \in \U_k(\HH)$ for every $k \in \N$. Furthermore, $\U_k(\HH) + \U_{\ell}(\HH) \subset \U_{k+\ell}(\HH)$, but the converse is not necessarily true. Moreover, $\ell \in \U_k(\HH)$ if and only if $k \in \U_\ell(\HH)$, and $1 \in \U_k(\HH)$ if and only if $k=1$, which is also equivalent to $\U_k(\HH) = \{1\}$.

In zero-sum theory over a finite group $G$, the monoid $\B(G)$ of product-one sequences is atomic, being Krull precisely when $G$ is abelian \cite[Proposition 3.4]{Oh}. In this case, $\B(G)$ is a natural model for studying the arithmetic of Krull monoids, and has been extensively investigated (see \cite{Ge,Sc2}). For non-abelian groups, $\B(G)$ is no longer Krull but remains a C-monoid \cite[Theorem 3.2]{CDG}, hence still enjoying finiteness properties for arithmetical invariants \cite{GeHK,GeZh} (see also \cite{Oh,Oh2}). 

For brevity, we write $*(\B(G))=*(G)$, where $* \in \{\rho, \rho_k, \lambda_k, \U_k, \LL, \A, \dots\}$. It is known that $\U_k(G)$ is the singleton $\{k\}$ if and only if $|G| \le 2$, and in this case we obtain that $\B(G)$ is half-factorial (see \cite[Proposition 3.3.2]{Ge}). Thus it is convenient to assume that $|G| \ge 3$. We have the following results.

\begin{teo}[{\cite[Theorem 5.5.1]{Oh}}]
	Let $G$ be a finite group with $|G| \ge 3$. Then for every $k \in \N$, $\U_k(G) = [\lambda_k(G),\rho_k(G)]$ is a finite interval.
\end{teo}

\begin{teo}[{\cite[Proposition 5.3]{OhZh1}}]
	Let $G$ be a finite group with $|G| \ge 3$. For every $\ell \in \N_0$, we have 
	$$\lambda_{\ell \, \D(G) + j}(G) = 
	\begin{cases}
		2\ell &\text{ for } j=0, \\
		2\ell+1 &\text{ for } j\in[1,\rho_{2\ell+1}(G)-\ell\, \D(G)], \\
		2\ell+2 &\text{ for } j\in[\rho_{2\ell+1}(G)-\ell \,\D(G)+1,\D(G)-1],
	\end{cases}$$
	provided that $\ell\,\D(G)+j > 0$.
\end{teo}

It is worth mentioning that if $G$ is infinite, then $\U_k(G) = \N_{\ge 2}$ \cite[Theorem 7.4.1]{GeHK}. In the context of the preceding theorems, $\rho_k(G)$ becomes a central invariant in the study of the interplay between zero-sum problems and factorization theory. In this direction, the following bounds hold.

\begin{prop}[{\cite[Lemma 1]{GeGrYu}}]\label{propdesigrho}
	Let $G$ be a finite group with $|G| \ge 3$.
	\begin{enumerate}[(i)]
		\item $k + \ell \le \rho_k(G) + \rho_\ell(G) \le \rho_{k+\ell}(G)$;
		\item $\rho_{2k}(G) = k\D(G)$ and 
		\begin{equation}\label{ineqkodd}
			k\D(G) + 1 \le \rho_{2k+1}(G) \le k\D(G) + \left\lfloor \frac{\D(G)}{2} \right\rfloor.
		\end{equation}
		In particular, $\rho(G)=\dfrac{\D(G)}{2}$.
		\item If $\rho_{2k+1}(G) \ge m$ for some $m \in \N$ and $\ell \ge k$, then $\rho_{2\ell+1}(G) \ge m + (\ell - k)\D(G)$.
	\end{enumerate}
\end{prop}

We observe that $\rho_k(G)$ is fully determined in terms of $\D(G)$ when $k$ is even. Nevertheless, for odd $k$, the lower bound on Inequality \eqref{ineqkodd} is attained for cyclic groups (see \cite[Corollary 1 and Proposition 6]{GeGrYu}), while the upper bound is conjectured to be eventually attained for non-cyclic abelian groups.

\begin{conj}[{\cite[Conjecture 1]{GeGrYu}}]
	Let $G$ be a finite non-cyclic abelian group with $\D(G) \ge 4$. Then there exists $k_0 \in \N$ such that $$\rho_{2k+1}(G) = k\D(G) + \left\lfloor \frac{\D(G)}{2} \right\rfloor$$
	for each $k \ge k_0$.
\end{conj}

By item (iii) of previous proposition, if this conjecture holds for some $k_0$, then it also holds for every $k \ge k_0$. Oh and Zhong investigated this problem for dihedral and dicyclic groups. In particular, they proved that the upper bound in Inequality \eqref{ineqkodd} is attained when $G$ is the dihedral group of order $2n$ with $n$ odd (see \cite[Theorem 5.4]{OhZh1}). On the other hand, for dihedral groups of order $2n$ with $n$ even, as well as for dicyclic groups of order $4m$, $m \ge 2$, they showed that for every $k \ge 2$, $\rho_k$ attains neither the lower nor the upper bound in Inequality \eqref{ineqkodd} (see \cite[Theorem 5.5]{OhZh1}).

For the group $C_q \rtimes C_p$, in this section we show that neither the lower nor the upper bound in Inequality \eqref{ineqkodd} is attained, a phenomenon similar to \cite[Theorem 5.5]{OhZh1}. This occurs because the extremal sequences described in Theorem \ref{main} somehow resemble those obtained for dihedral groups of order $2n$ with $n$ even and for dicyclic groups. The main result of this section is stated as follows.

\begin{teo}
	Let $p,q$ be odd prime numbers with $p \mid q-1$ and let $\G \cong C_q \rtimes C_p$. For every $k\in\N$, we have that 
	$$k\D(\G)+2\le\rho_{2k+1}(\G)\le k\D(\G)+\dfrac{\D(\G)}{2}-1.$$ 
\end{teo}

\proof 
	Recall that $\D(\G) = 2q$. From Proposition \ref{propdesigrho}(iii), in order to prove the first inequality it suffices to show that $\rho_3(\G) \ge 2q+2$. As a consequence of Theorem \ref{inv}, we consider the minimal product-one sequences 
	$$S_1 = y^{[2q-2]}\bd x\bd x^{-1}y^{s^{p-1}+1}, \; S_2 = (y^{-1})^{[2q-2]}\bd x^{-1}y^{-1}\bd xy^{-1} \in \A(\G).$$ 
	Moreover, since the products
	$$\begin{aligned}
		x^{-1} \cdot xy^{-s-1} = y^{-s-1}, && x^{-1} \cdot xy^s = y^s, && x^{-1}y^{s^{p-1}} \cdot xy^{-s-1} = y^{-s}, && x^{-1}y^{s^{p-1}} \cdot xy^s = y^{s+1}
	\end{aligned}$$
	are all different from $1$, it follows that 
	$$S_3 = x^{-1}\bd xy^{-s-1}\bd xy^s\bd x^{-1}y^{s^{p-1}} \in \A(\G)$$
	is a minimal product-one sequence as well. We obtain a distinct factorization
	$$S_1\bd S_2\bd S_3 = U_1\bd U_2\bd U_3\bd U_4\bd U_5^{[2q-2]},$$ 
	where 
	$$U_1 = x\bd x^{-1}, \; U_2 = x^{-1}y^{s^{p-1}+1}\bd xy^{-s-1}, \; U_3 = x^{-1}y^{-1}\bd xy^s, \; U_4 = x^{-1}y^{s^{p-1}}\bd xy^{-1}, \; U_5 = y\bd y^{-1} \in \A(\G)$$
	are minimal product-one sequences. This implies that $2q+2 \in \mathcal{U}_3(\G)$, whence $\rho_3(\G) = \sup \U_3(\G) \ge 2q+2$.

	For the upper bound, we assume that $\rho = \rho_{2k+1}(\G) = q(2k+1)$ for some $k \in \N$. Suppose in addition that $k$ is minimal with this property. By assumption, there exist minimal product-one sequences $V_1, \dots, V_{2k+1} \in \A(G)$ 
	such that 
	$$\rho \in \L(V_1 \bd \ldots \bd V_{2k+1}).$$ 
	By definition, there exist minimal product-one sequences $W_1, \dots, W_\rho \in \A(G)$ 
	such that 
	$$T = V_1 \bd \ldots \bd V_{2k+1} = W_1 \bd \ldots \bd W_\rho.$$

	If $1^{[2]} \mid T$, then the sequence $T \bd (1^{[2]})^{[-1]}$ contradicts the minimality of $k$. If $1 \mid T$ but $1^{[2]} \nmid T$, say $V_{2k+1} = W_{\rho} = 1$, then we consider $T \bd 1^{[-1]} = V_1 \bd \dots \bd V_{2k} = W_1 \bd \dots \bd W_{\rho-1}$. Since $\D(\G) = 2q$ and $1^{[2]} \nmid T$, it follows that $|V_i| \le 2q$ and $|W_j| \ge 2$ for every $i \in [1,2k]$ and $j \in [1,\rho-1]$. Thus
	$$4qk \ge |V_1 \bd \dots \bd V_{2k}| = |W_1 \bd \dots \bd W_{\rho-1}| \ge 2(\rho-1) = 4qk+2q-2,$$
	a contradiction since $q \ge 7$. Hence $1 \nmid T$ and this implies that $|W_j| \ge 2$ for every $j \in [1,\rho]$. Since $|V_i|\le 2q$ for every $i \in [1,2k+1]$, it follows that 
	$$2q(2k+1) \ge |V_1 \bd \ldots \bd V_{2k+1}| = |W_1 \bd \ldots \bd  W_\rho| \ge 2q(2k+1),$$ 
	whence $|V_i|=2q = \D(\G)$ and $|W_j|=2$ for every $i \in [1,2k+1]$ and $j \in [1,\rho]$. By Theorem \ref{main}, each $V_i$ has an unique element $g_i \in \G$ with $\ord(g_i) = q$ such that $g_i \in \supp(V_i)$. More precisely, $g_i^{[2q-2]} \mid V_i$. On the other hand, since $|W_j|=2$, if $g^{[2q-2]}\mid V_i$ for some $i \in [1,2k+1]$, then there exists $\ell$ such that $(g^{-1})^{[2q-2]} \mid V_\ell$. Since $2k+1$ is odd, there exists $V_i$ such that none of the terms of order $q$ can be paired with its inverse to form $W_j$ for some $j$. This leads to the desired contradiction, thereby completing the proof of the theorem.

\qed

\end{document}